\documentclass[12pt]{article}
\usepackage{amsmath,amsthm,amsfonts}

\input epsf.tex
\newdimen\epsfxsize
\newdimen\epsfysize
\def\qed{\vrule height5pt width3pt depth.5pt}

\theoremstyle{plain}
\newtheorem{thm}{Theorem}[section]

\newtheorem{lem}[thm]{Lemma}
\newtheorem{prop}[thm]{Proposition}

\newtheorem{conj}{Conjecture}[section]
\newtheorem{rem}{Remark}[section]

\begin{document}

\title{Virtual knots undetected by 1 and 2-strand bracket
polynomials }

% infomation for author
\author{H. A.  Dye \\
  MADN-MATH \\
  United States Military Academy \\
  646 Swift Road \\
  West Point, NY 10996-1905 \\
 hdye@ttocs.org }

\maketitle

\begin{abstract}
Kishino's knot  is not detected by the
 fundamental 
group or the bracket polynomial; these invariants cannot differentiate between Kishino's knot and the unknot. However, we can show that Kishino's knot is not equivalent to unknot 
by applying either the 3-strand bracket polynomial or the surface bracket polynomial. In this paper, we construct two non-trivial virtual knot diagrams, $ K_D $ and $ K_m $, that are not 
not detected by the bracket polynomial or the 2-strand bracket polynomial. From these diagrams, 
we construct two infinite families of non-classical virtual knot diagrams that are not detected by the bracket polynomial. Additionally, we note these virtual knot diagrams are trivial as flats.
\end{abstract}

\section{Introduction}
\begin{figure}[htb] \epsfysize = 0.75 in
\centerline{\epsffile{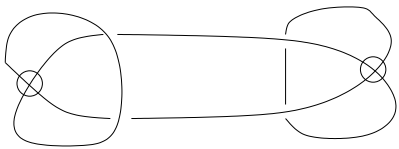}}
\caption{Kishino's Knot}
\label{fig:kishknot}
\end{figure}

   Kishino's knot, illustrated in Figure \ref{fig:kishknot}, is not
detected by the fundamental group, bracket polynomial or the 2-strand 
bracket polynomial. Kishino and Satoh  
\cite{kishpoly} demonstrated that Kishino's knot is detected by the 3-strand bracket 
polynomial and that this virtual knot diagram is not equivalent to the unknot. 
The surface bracket polynomial \cite{dk1} indicates that this knot is non-classical - not equivalent to a classical knot diagram. We construct other examples of this phenomena in this paper.

We review virtual knot theory and recall the definition of the $ N$-strand bracket polynomial, the fundamental group of a virtual knot diagram, and the surface bracket polynomial \cite{dk1}.  
The fundamental group and bracket polynomial are invariants that are commonly used to determine if a classical knot diagram is equivalent to the unknot. The surface bracket polynomial can be applied to show that virtual knot diagrams are non-classical.

We introduce the virtual knot diagram, $ K_D $, which is not detected by the bracket polynomial or
the 2-strand bracket polynomial. 
The 3-strand bracket polynomial shows that $ K_D $ is non-trivial and the surface bracket polynomial demonstrates that this knot is non-classical.  This knot diagram is used to construct
an infinite family of non-classical virtual knot diagrams that are not detected the bracket polynomial (Family A). We apply the surface bracket polynomial to show that all members of the family are non-classical. 

We modify $ K_D $ and Kishino's knot to construct the virtual knot diagram $ K_m $. This knot diagram is not detected by the fundamental group, bracket polynomial, or 2-strand bracket polynomial. The 3-strand bracket polynomial detects this virtual knot diagram. From $ K_m $, we construct an infinite family of non-classical virtual knot diagrams that are not detected by the fundamental group or the bracket polynomial (Family B). 

  $ K_D $ and $ K_m $ are both detected by  
the 3-strand bracket polynomial. Computing the 3-strand bracket polynomial is a labor intensive process  (motivating a search for new invariants such as the surface bracket polynomial \cite{dk1}). We conjecture that the 3-strand bracket polynomial detects all the virtual knot diagrams in these families. Note that if the bracket polynomial or the 2-strand bracket polynomial detects a virtual knot diagram then the 3-strand bracket polynomial will also detect this diagram. As a result, we conjecture that the 3-strand bracket polynomial will detect all virtual knot diagrams.

\section{Virtual Knot Diagrams} 

\begin{figure}[htb] \epsfysize = 1 in
\centerline{\epsffile{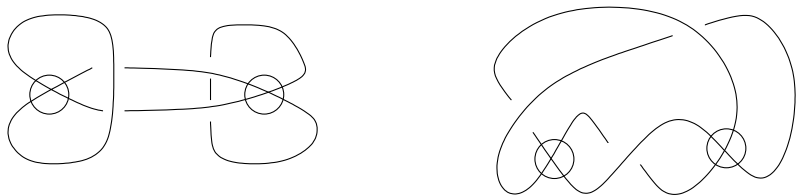}}
\caption{Examples of Virtual Knot Diagrams}
\label{fig:virtex}
\end{figure}

A \emph{virtual knot diagram} is a decorated immersion of $ S^1 $ in the 
plane. A virtual knot diagram has two types of crossings: classical crossings and virtual crossings.
We
indicate classical crossings with over-under markings and the virtual crossings are indicated by 
a solid encircled X. Two virtual knot diagrams are illustrated in Figure \ref{fig:virtex}.
Note that the classical knot diagrams are a subset of the 
virtual knot diagrams. 

\begin{figure}[htb] \epsfysize = 0.75 in
\centerline{\epsffile{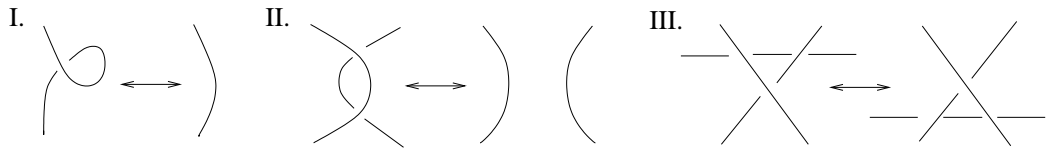}}
\caption{Reidemeister Moves}
\label{fig:rmoves}
\end{figure}

\begin{figure}[htb] \epsfysize = 1.75 in
\centerline{\epsffile{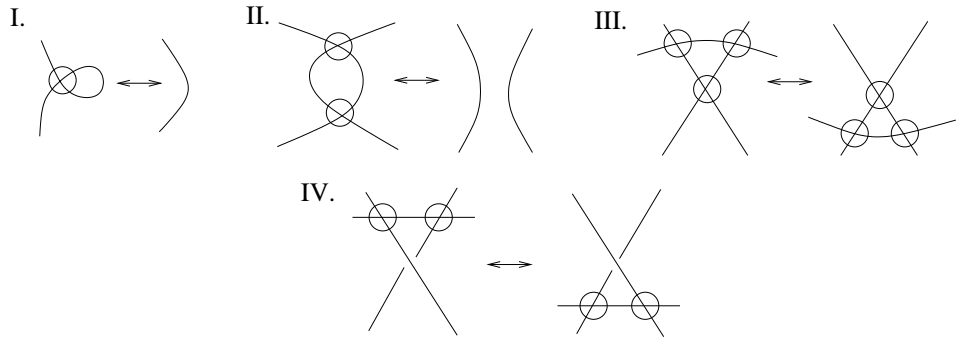}}
\caption{Virtual Reidemeister Moves}
\label{fig:vrmoves}
\end{figure}

We recall the \emph{Reidemeister moves}. Local versions of the classical Reidemeister moves are illustrated in Figure \ref{fig:rmoves}.
Two classical knot diagrams are said to be \emph{equivalent} if one may be transformed
into the other by a sequence of Reidemeister moves. 
To extend the notion of equivalence to virtual knot diagrams, we extend 
our set of diagrammatic moves to 
include virtual crossings. The \emph{virtual Reidemeister moves} 
are illustrated in Figure \ref{fig:vrmoves}.
Note that the virtual Reidemeister moves I, II, and III involve only virtual crossings. 
Two virtual knot diagrams are said to be \emph{virtually equivalent} if one diagram
may be transformed into the other via a sequence of Reidemeister and virtual Reidemeister moves. 

\begin{figure}[htb] \epsfysize = 0.5 in
\centerline{\epsffile{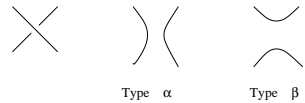}}
\caption{Smoothing Types}
\label{fig:smooth}
\end{figure}

To introduce the  generalized bracket polynomial of a virtual knot diagram, we define a smoothing of a classical crossing and a state of a virtual knot diagram. 
We \emph{smooth} a classical crossing in a virtual knot diagram by removing a small neighborhood of the classical crossing and replacing it with neighborhood containing two non-intersecting segments. The classical crossing is replaced with either a type $ \alpha $  
smoothing or a type $ \beta $ smoothing as shown in Figure \ref{fig:smooth}.
(Each smoothed crossing is implicitly labeled with it's 
smoothing type.) We obtain a
\emph{state} of a virtual knot diagram by assigning a smoothing type to each classical crossing and smoothing the diagram accordingly.
A state of a virtual knot diagram consists of a set of closed curves that possibly contain 
virtual crossings. The
set of all states of a virtual knot diagram is denoted as $ S $.
 Note that a virtual knot diagram with 
$ N $ classical crossings and $ M $ virtual crossings
has $ 2^N $ states.

We define the \emph{generalized bracket polynomial} of a virtual knot 
diagram.
Let $ d = -A^{-2} -A^2 $. For a given state $ s \in S $, let $ c(s) $ equal the 
number of type $ \alpha $ smoothings minus the number of type 
$ \beta $ smoothings. Let $ | s | $ represent the number of closed
curves in the state, $ s $. 
We
denote the generalized bracket polynomial of $ K $ as $ \langle K \rangle $, 
then:
\begin{equation}
\langle K \rangle = \underset{s \in S}{ \sum } A^{ c(s)} d^{ |s| -1 }
\end{equation}
This polynomial is invariant under the Reidemeister II and III moves and the 
virtual Reidemeister moves \cite{kvirt}. We will refer to this polynomial as 
the bracket polynomial for the remainder of the paper.

\begin{figure}[htb] \epsfysize = 2 in
\centerline{\epsffile{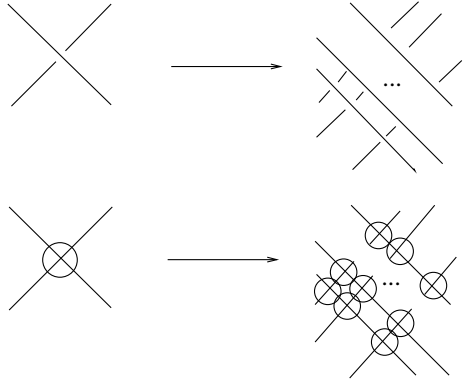}}
\caption{Relationship between $ K $ and 
$  K_N $ at crossings}
\label{fig:Nstrand}
\end{figure}

Let $ K $ be a virtual knot diagram. We modify $ K $ to form the virtual link diagram $ K_N $ for $ N \geq 1 $. 
The diagram $ K_N $ is formed by taking $ N $ parallel copies of 
$ K $. The 
relationship between these parallel copies at the classical and virtual crossings is illustrated in Figure \ref{fig:Nstrand}.

Note that the virtual knot diagram $ K_1 $ is  $ K $. 
We use $ K_N $ to define the \emph{N-strand bracket polynomial} of $ K $. 
We denote the N-strand bracket 
polynomial of a virtual knot diagram $ K $  as $ \langle K \rangle_N $ then:
\begin{equation}
 \langle K \rangle_N = \langle K_N \rangle. 
\end{equation}

\begin{figure}[htb] \epsfysize = 0.5 in
\centerline{\epsffile{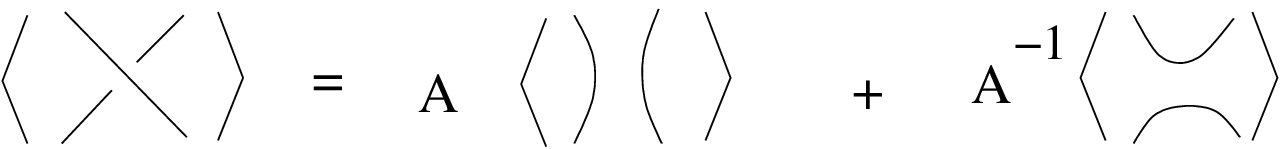}}
\caption{Skein Relation}
\label{fig:skein}
\end{figure}

The bracket polynomial and the N-strand bracket polynomial may also be computed using the skein relation shown in Figure \ref{fig:skein}.

The N-strand bracket polynomial is invariant under the Reidemeister II and III 
moves and the virtual Reidemeister moves. In particular, any move performed on 
the diagram $ K $ can be replicated on the diagram $ K_N $ by a sequence of
the same move. The N-strand bracket polynomial is not invariant under the Reidemeister I move.

\begin{figure}[htb] \epsfysize = 1 in
\centerline{\epsffile{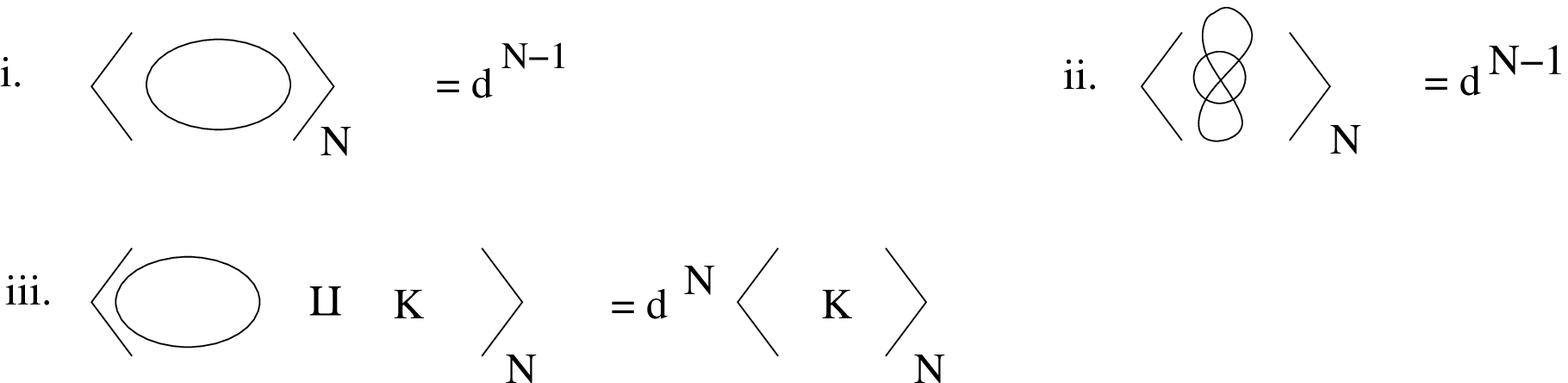}}
\caption{Evaluation of N-strand bracket polynomial }
\label{fig:genstrand}
\end{figure}

Several key facts about the N-strand bracket polynomial are illustrated in Figure \ref{fig:genstrand}.

Note that as $ N $ increases the computational complexity of this invariant 
grows exponentially. 
A diagram with $ n $ classical crossings has $ 2^n $ states, but 
$ K_2 $ has $ 2^{4n} $ and $ K_3 $ has $ 2^{9n} $ states.
We will restrict our attention to the 2-strand and 3-strand bracket polynomial.
For classical knot diagrams, the following conjecture has been made about the bracket polynomial.
\begin{conj}For a classical knot diagram $ K $, if $ \langle K \rangle = 1 $ then $ K $ is the unknot. 
\end{conj}
However, for virtual knot diagrams the conjecture is false. Kauffman \cite{kvirt} demonstrated that there are an infinite number of virtual knot diagrams, $ K $, such that
$ \langle K \rangle = 1 $. 

\begin{figure}[htb] \epsfysize = 1 in
\centerline{\epsffile{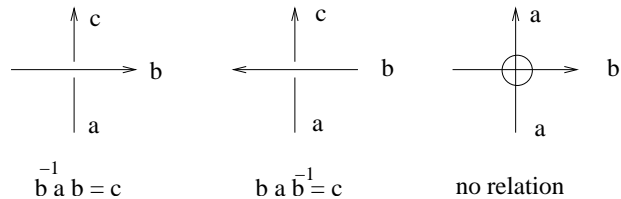}}
\caption{Fundamental Group Relations}
\label{fig:fund}
\end{figure}

The \emph{fundamental group} \cite{kvirt} of a virtual knot diagram  is computed from
a labeled oriented diagram. (See \cite{purp} for the classical definition of fundamental group.)
Let $ K $ be an oriented virtual knot diagram with $ n $ classical crossings  and $ 2n $ arcs. The \emph{arcs} in a virtual knot diagram have
endpoints at the classical crossings and pass through virtual 
crossings without termination.
The fundamental group of
$ K$, denoted $ \pi_1 (K ) $, 
 is the free group generated by the labels on the arcs modulo relations 
determined by the classical crossings in diagram. Each crossing produces one of the relations 
illustrated in Figure \ref{fig:fund}.
For example, a knot with $ n$ classical crossings will have $ 2n $ generators and $ n $ relations.

Kishino's knot is not detected by the fundamental group. The fundamental group of Kishino's knot is $ \mathbb{Z} $ which is equivalent to  the fundamental group of the unknot.

We recall the surface bracket polynomial from \cite{dk1}.
Referring to \cite{dk1}, \cite{kamada}, and \cite{kvirt}, we recall that equivalence classes of virtual knot diagrams are in one to one correspondence with equivalence classes of knot diagrams on two dimensional surfaces. 
Two surface-knot diagram pairs are equivalent when one may be transformed into the other via a sequence of Reidemeister moves in the surface, homeomorphisms of the surfaces, and handle cancellations or additions. 
A fixed oriented 2-dimensional surface $ F $ with an immersed knot diagram  $ K $ is a \emph{representation} of a virtual knot diagram. From a representation, we recover the virtual knot diagram $ \hat{K} $ by projecting the diagram onto  the plane. Virtual crossings result when a double point in the projection does not correspond to a classical crossing in the immersed diagram. We construct a representation of a virtual knot diagram by the following process.
The virtual knot diagram may be viewed as a decorated immersion
 on the surface of a sphere instead of a plane.
For each virtual crossing in the diagram, remove a small neighborhood of one of the arc passing through the crossing. We attach a handle with an appropriate embedding of an arc to the sphere.
 
We denote a specific representation with surface $ F $ and embedding $ K $ as $ (F, K ) $.

A representation of the virtual knot diagram $ \hat{K} $ is said to have \emph{minimal genus} if the surface $ F $ has the minimum genus of all representations $ (F, K) $ of $ \hat{K} $. 

We realize the following lemma from \cite{kup}.
\begin{lem}Let $ (F,K) $ be a representation of $ \hat{K} $. If the minimal genus of $ \hat{K} $ is greater than zero than $ \hat{K} $ is non-trivial and non-classical. 
\end{lem}
This lemma allows us to use the surface bracket polynomial to determine if a virtual knot diagram is non-trivial and non-classical in some cases. We introduce the surface bracket polynomial with 
the definition of the \emph{surface-state} pair of a representation.

From $ (F,K) $, we obtain the surface-state pair $ (F,s ) $ by choosing a smoothing type for each classical crossing.
This results in a set of simple closed curves (possibly bounding a disk) of the surface $ F $

 We define the \emph{surface bracket polynomial} of a 
representation  $ (F,K) $ \cite{dk1}.
 Let $ \hat{K} $ be a
 virtual knot diagram,
and let  $ (F, K ) $ be a fixed representation of $ 
\hat{K} $. 
We denote the surface bracket polynomial of $ K $ 
 as $ \langle (F,K) \rangle $. Then:
$$ 
 \langle (F,K) \rangle = \underset{(F,s(c)) \in (F,S)}{ \sum} \langle K| s(c) \rangle d^{|s(c)|} [s(c)]
$$
where $ \langle K | s(c) \rangle = A^{c(s)} $ and $ c(s) $ is the number of type $ \alpha
$ smoothings minus the number of type $ \beta $ smoothings. 
$ | s(c) | $ is the number of curves which bound a disk in the surface 
and $ [ s(c) ] $ represents a formal sum
 of the disjoint curves that do not bound a disk in 
the surface-state pair $ (F,s(c)) $.
We apply the following theorem \cite{dk1} to determine if a virtual knot diagram is non-trivial and non-classical using the surface bracket polynomial. 

\begin{thm}Let $ (F, K ) $ be a representation of a virtual knot diagram with 
$ F = T_1 \sharp T_2 \ldots 
\sharp T_n $. 
Let 
$$ \lbrace ( F, s_1), (F, s_2) \ldots (F,s_m) \rbrace $$ 
denote  the 
collection of surface-state pairs obtained from $ (F,K) $. 
Assign an arbitrary orientation to each curve in the surface-state
pairs.
Let $p: F \rightarrow T_k $ be the collapsing map, and let
$ p_*: H_1 (F, \mathbb{Z}) \rightarrow H_1 (T_k, \mathbb{Z} ) $
be the induced map on homology.  
If for each $ T_k $ there exist two states $ s_i $ and 
$ s_j $  with non-zero coefficients that contain curves 
(with arbitrarily assigned orientation) $ \gamma_i $ and $ \gamma_j $ respectively, 
such that $ p_*[ \gamma_i] \bullet p_*[ \gamma_j]  \neq 0 $ then
there is no cancellation curve for $ (F,K) $.
\end{thm}
\begin{rem}As a result of this theorem, if no cancellation curve exists on the surface of the representation then this is a minimal genus representation.
\end{rem}

In this paper, we demonstrate that are an infinite number of  non-classical virtual knot diagrams that are not detected by the bracket polynomial and the fundamental group. We will construct two families of virtual knot diagrams. The first, Family A, is not detected by the bracket polynomial. 
The second, Family B, is a modification of Family A and Kishino's knot. 
This family is not detected by the bracket polynomial or the fundamental group.
 The surface bracket polynomial shows that both these families are non-classical.
\section{Infinite Family A}

Kishino's knot, illustrated in Figure \ref{fig:kishknot}, is not detected by the bracket polynomial or the 2-strand bracket polynomial. This knot \cite{kishpoly} was determined to have
a non-trivial 3-strand bracket polynomial, proving that Kishino's knot is not equivalent to the unknot. In this section, we construct a non-trivial virtual knot diagram, $ K_D $, that is not detected by the bracket polynomial or the 2-strand bracket polynomial. Using $ K_D $, we construct an infinite family (Family A) of non-trivial virtual knot diagrams that are not detected by the bracket polynomial. 
We prove that this family of diagrams is non-trivial using the surface bracket polynomial \cite{dk1}.
 We conjecture that family A is not detected by the N-strand bracket polynomial when $ N < 3 $. This conjecture has been verified for the first element of this family. We are unable to prove this conjecture for all members of this family due to the complexity of computing the 3-strand bracket polynomial of diagrams with even a few classical crossings. 
Note that the 3-strand bracket polynomial of  Kishino's knot  has $ 2^{36} $ states. The number of states increase exponentially in proportion to the number of classical crossings. For the 3-strand bracket polynomial, each additional classical crossing increases the number of states by the multiple $ 2^9 $. 
\begin{rem} Kishino's knot has been detected by the quaternionic biquandle \cite{bart}, Kadokami's methods \cite{kadokami}, and the surface bracket polynomial \cite{dk1}. Kishino's knot is non-trivial as a flat diagram on a surface and, as a result, is detected by Kadokami's methods. 
\end{rem}

\begin{figure}[htb] \epsfysize = 1.5 in
\centerline{\epsffile{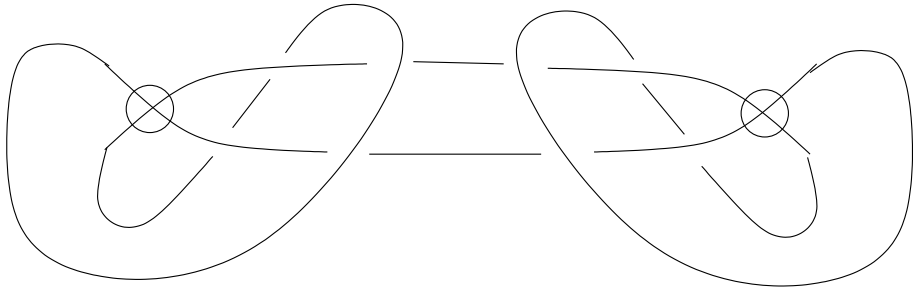}}
\caption{Virtual Knot Diagram, $ K_D $}
\label{fig:heatherknot}
\end{figure}

\begin{prop}
The virtual knot diagram, $ K_D $, illustrated in  Figure \ref{fig:heatherknot} is non-trivial. This virtual knot diagram is not detected by the bracket polynomial or the 2-strand bracket polynomial. However, the 3-strand bracket polynomial detects $ K_D $ since $  \langle K_D \rangle_3 \neq d^2 $.
\end{prop}

\textbf{Proof:}
Note that:
\begin{align*}
\langle K_D \rangle &= 1 & \qquad \langle K_D \rangle_2 &= d
\end{align*}
However:
% Shift description to A terms
\begin{gather*}
\langle K \rangle_3  = d (
-1842 + ( A^{76} + A^{-76} ) -2(A^{ 72}+A^{-72})\\ 
-2(A^{70} + A^{-70})-3(A^{68} + A^{-68} ) -8(A^{66} + A^{-66}) \\
 -18(A^{64} + A^{-64})-17(A^{62}+ A^{-62}) -8 ( A^{60} + A^{-60}) \\
-5(A^{58} + A^{-58})-2(A^{56} + A^{-56}) + 21(A^{54} + A^{-54}) \\ 
+ 66(A^{52} + A^{-52} )  + 95(A^{50} + A^{-50}) + 103(A^{48} + A^{-48}) \\ 
 +114(A^{46} + A^{-46}) +136 ( A^{44} + A^{-44} ) + 111( A^{42} + A^{-42} ) \\
 + 11( A^{40} + A^{-40})-110(A^{38} + A^{-38} ) -209(A^{36} + A^{-36} )\\
 -326(A^{34} + A^{-34} ) -491(A^{32} + A^{-32}) -601(A^{30} + A^{-30}) \\
-559(A^{28} + A^{-28}) -380(A^{26} + A^{-26}) -178(A^{24} + A^{-24}) \\
+142(A^{22} + A^{-22})  +594(A^{20} + A^{-20}) +1050(A^{18} + A^{-18})\\
 +1329(A^{16} + A^{-16})+1334(A^{14}+A^{-14}) + 1215(A^{12} + A^{-12}) \\
+ 814(A^{10} + A^{-10} ) + 193(A^{8} + A^{-8}) -573(A^{6} + A^{-6})\\
 -1257(A^{4} + A^{-4}) -1660( A^{2} + A^{-2})) 
\end{gather*}
Hence, $ K_D $ is non-trivial.
\qed

\begin{figure}[htb] \epsfysize = 2.5 in
\centerline{\epsffile{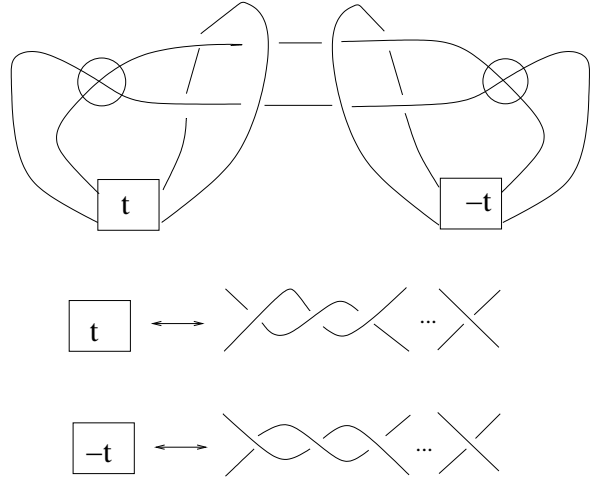}}
\caption{Schematic of Infinite Family A, $ K_D(t) $ }
\label{fig:schem}
\end{figure}

We introduce an infinite family of virtual knot diagrams based on this diagram (Family A). These diagrams in Family A are denoted by $ K_D (t) $ where $ t $ $ (-t) $ represents the number of inserted positive (negative) twists, as shown in Figure \ref{fig:schem}.
 The members of Family A are not detected by the bracket polynomial.

\begin{figure}[htb] \epsfysize = 2.5 in
\centerline{\epsffile{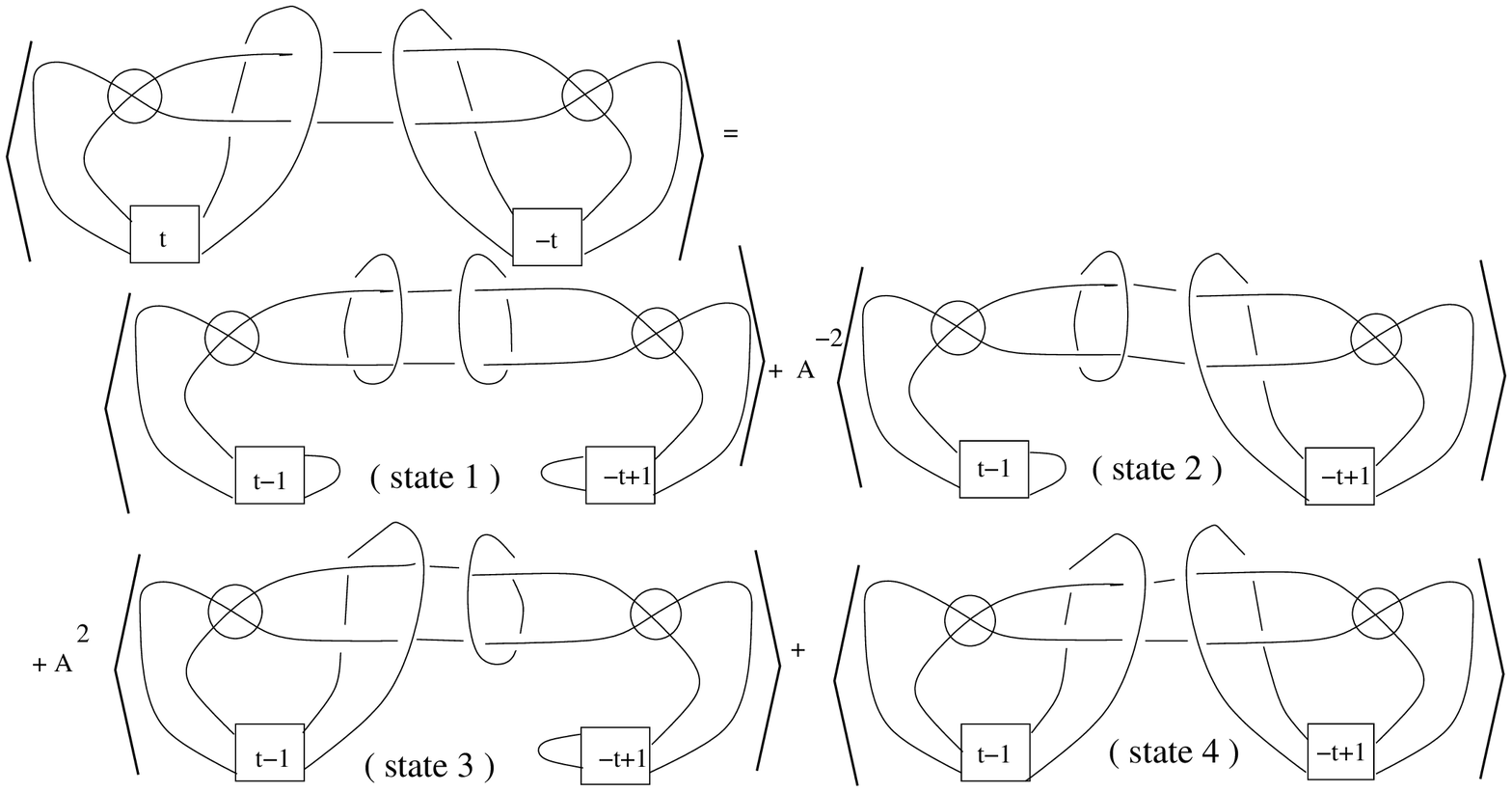}}
\caption{Skein Relation Applied to Family A}
\label{fig:sexpan}
\end{figure}

\begin{thm} \label{bkd} The virtual knot diagrams, $ K_D (t) $, in Family A as illustrated in Figure \ref{fig:schem} are not 
detected by the bracket polynomial. 
\end{thm}
\textbf{Proof:}Let $ K_D(t) $ represent the virtual knot diagram with  $ t $ twists on the left and $ -t $ twists inserted on the right.
Note that if $ t=0 $, $K_D(0)$ is the diagram $ K_D $. 
We recall that $ \langle K_D(0)\rangle = 1 $. 
Expanding the diagram $ K_D(t)$ using the skein relation, we obtain  the sum shown in Figure \ref{fig:sexpan}.

We assume that $  \langle K_D(t-1) \rangle = 1 $.
Note that the first diagram (state 1) in the expansion is equivalent to three 
unlinked loops after an appropriate sequence of virtual and Reidemeister moves. This sequence of moves does not change the writhe of the diagram. As a result, the bracket polynomial of this diagram is $ d^2 $. Similarly,  state 2 and state 3 are equivalent to two unlinked loops and the bracket polynomial of these diagrams is $ d $. We obtain:
\begin{gather*}
\langle K_D (t) \rangle = \langle \text{state 1} \rangle + 
A^{-2} \langle \text{state 2} \rangle + A^2\langle \text{state 3} \rangle + \langle K_D (t-1) \rangle  \\
\langle K_D(t) \rangle = d^2 + A^{-2} d + A^{2} d + 1
\end{gather*}
This reduces to:
\begin{gather*}
\langle K_D(t) \rangle = A^4 + 2 + A^{-4} - 1 - A^{-4} - A^4 - 1 + 1 \\
\langle K_D(t) \rangle = 1  
\end{gather*}
Hence, no member of Family A is detected by the bracket polynomial.
\qed

\begin{conj} We conjecture that the family depicted in Figure \ref{fig:schem} is detected by the 3-strand bracket polynomial
but not the bracket polynomial or the 2-strand bracket polynomial.
\end{conj}

These computations have been verified for the virtual knot diagrams $ K(0) $ and $ K(1) $. The 3-strand bracket polynomial only determines that these diagrams are non-trivial and does not show that they are non-classical. We are unable to verify that the 3-strand bracket polynomial detects these virtual knot diagrams when
 when $ t \geq 2 $ due to the large number of computations involved.

We apply the surface bracket polynomial to show that the members of Family A are not equivalent to a classical knot diagram.

\begin{figure}[htb] \epsfysize = 1 in
\centerline{\epsffile{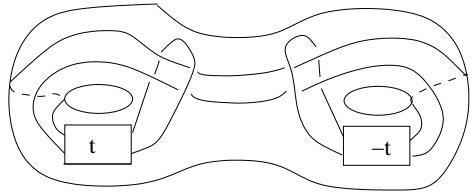}}
\caption{Schematic Representation of $ K_D (t)$ }
\label{fig:kdsurf}
\end{figure}

\begin{thm} All the virtual knot diagrams in Family A, shown in Figure \ref{fig:schem}, are non-trivial and non-classical.
\end{thm}
\textbf{Proof:}
In Figure \ref{fig:kdsurf}, we illustrate a schematic representation of $ K_D (t) $ in the connected sum of two tori, F. 

\begin{figure}[htb] \epsfysize = 1 in
\centerline{\epsffile{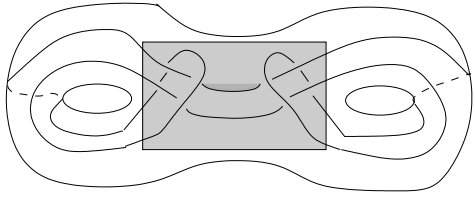}}
\caption{A Representation of $ K_D (0)$ }
\label{fig:kdsurf0}
\end{figure}

We compute the surface bracket polynomial of  $ (F, K_D (0)) $ based on the representation shown in Figure \ref{fig:kdsurf0}.

\begin{figure}[htb] \epsfysize = 2.5 in
\centerline{\epsffile{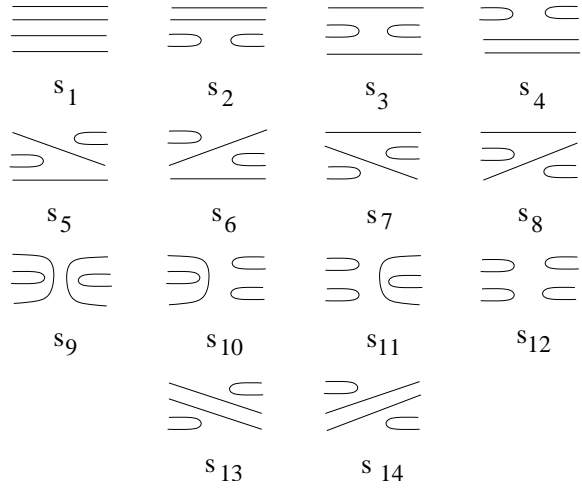}}
\caption{Possible bracket expansions of a 4-4 tangle}
\label{fig:expan4}
\end{figure}

To compute the surface bracket polynomial, we expand the 4-4 tangle contained in the shaded box in Figure \ref{fig:kdsurf0}. 
Note that expanding a classical $ 4-4 $ tangle via the skein relation results in one of the 14 states shown in Figure \ref{fig:expan4}. (These states are all elements of the $ 4^{th} $ Temperly-Lieb algebra.)

\begin{figure}[htb] \epsfysize = 1 in
\centerline{\epsffile{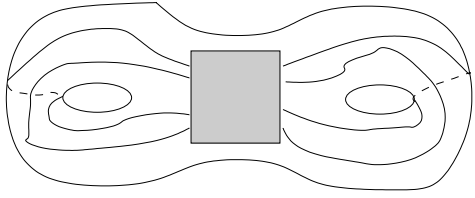}}
\caption{Expanded states of $K_D (t)$}
\label{fig:kdexpand}
\end{figure}

The tangle in the shaded box 
expands into 13 states with coefficients in 
$ \mathbb{Z}[A, A^{-1}] $. These states are placed in shaded
box of Figure \ref{fig:kdexpand} to form 13 surface-state pairs with non-zero coefficients.

Using this expansion, we determine that:
\begin{gather*} 
\langle (F,K_D (0) ) \rangle =
(-1 ) \langle (F, s_2) \rangle
+ ( A^{-8}-2 + A^8 ) \langle (F, s_3) \rangle
+(A^{-4} -2 + A^4) \langle (F, s_4) \rangle \\
+(A^{-4}-A^{-2} -A^2 + A^6 ) \langle (F, s_5) \rangle
+(A^{-6}-A^{-2}-A^2+A^6 ) \langle (F, s_6) \rangle \\ 
+ (A^{-10} + A^{-2}) \langle (F, s_7) \rangle 
+(-A^2 +A^{10}) \langle (F, s_8 ) \rangle 
+ (A^{-6}-A^{-2}-A^2 + A^6) \langle (F, s_9) \rangle \\
+ (A^{-4} -2+ A^8) \langle (F, s_{10}) \rangle 
+ (A^{-8} -2+A^4 ) \langle (F, s_{11}) \rangle
+( A^{-6}-A^{-2} -A^2 + A^6) \langle (F, s_{12}) \rangle \\
+ ( A^{-8} -A^{-4} ) \langle (F, s_{13}) \rangle 
+ ( -A^{4} + A^8 ) \langle (F, s_{14}) \rangle
\end{gather*}
The states $ s_3, s_4, $ and $ s_5 $ have non-zero coefficients. The presence of these states is sufficient to prevent the existence of a cancellation curve. This indicates that the  
 minimal genus of $ K_D (0) $ is two.
To compute the surface bracket polynomial of $ K_D (t) $ for $ t \geq 1 $, we use the expansion given in the proof of Theorem \ref{bkd}. We obtain a sum of representations of the virtual knots and links shown in Figure  \ref{fig:kdexpand}. One of these surface-link pairs is a 
 representation of $ K_D (t-1) $. Repeated applications of the skein 
relation  result in a representation of $ K_D (0 ) $.  The expansion of the other surface-link pairs does not result in  states
 $ (F, s_3) , (F, s_4) $ or $ (F, s_5) $.  If $ K_D (t) $ is classical, the coefficients of these states in the final expansion must be zero in order to admit a cancellation curve. However, the existence of a representation of $ K_D (0) $ as a partially expanded surface-state indicates that these states have non-zero coefficients. 
 Therefore, $ K_D (t) $ is non-trivial and non-classical. 
\qed

\begin{thm}
The virtual knot diagram in Figure \ref{fig:heatherknot} has a non-trivial fundamental group.
\end{thm}
\textbf{Proof:}
Computation determines:
\begin{equation*}
\pi_1 (K_D) = \lbrace a,g | a^{-1}ga = g^{-1}ag \rbrace
\end{equation*}
We  define  $ \rho : \pi_1 (K_D ) \rightarrow  GL_3 ( \mathbb{Z} ) $ 
to demonstrate that this group is non-trivial.
\begin{align*}
\rho (a) = \begin{bmatrix} 0 & 1 & 0 \\
			1 & 0 & 0 \\
			0 & 0 & 1 \end{bmatrix}  
 \qquad
\rho (g) =\begin{bmatrix} 0 & 0 & 1 \\
			0 & 1 & 0 \\
			1 & 0 & 0 \end{bmatrix}  
\end{align*}
Some additional computations will demonstrate that this is a non-Abelian representation.
We compute a general formula for $ \pi_1 (K_D(t)) $.
\begin{gather*}
x_0 = g  \text{ and } y_0  = a \\
\text{For t odd: } x_t  = y_{t-1}  \text{ and } y_t  = (y_{t-1})^{-1} x_{t-1} y_{t-1}  \\
\text{For t even: } x_t = y_{t-1} \qquad y_t  = y_{t-1} x_{t-1} (y_{t-1})^{-1} \\
\end{gather*}
Using this relation, we determine that:
\begin{equation*}
\pi_1 (K_D (t)) = \lbrace a^{-1} x_t a = g^{-1} y_t g \rbrace
\end{equation*}
We apply the mapping defined above to $ \pi_1 (K_D (t)) $. We note that $ \pi_1 (K_D (t)) $ maps into a non-trivial, non-Abelian group.
\qed

As a result of this theorem, we note that all virtual knot diagrams in Family A have non-trivial fundamental groups unlike Kishino's knot. 

\begin{rem}Although Family A is detected by the fundamental group, each diagram in Family A is trivial as a flat virtual knot diagram and is not detected by 
Kadokami's methods \cite{kadokami}.
\end{rem}  

In the next section, we construct a family that is not detected by the fundamental group or the bracket polynomial. 

\section{Infinite Family B}
\begin{figure}[htb] \epsfysize = 1.5 in
\centerline{\epsffile{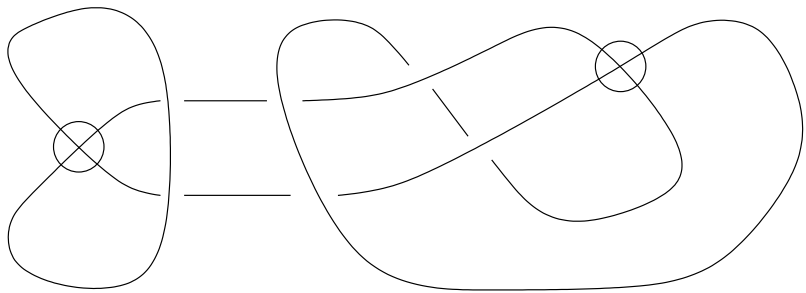}}
\caption{Virtual Knot Diagram, $ K_m $}
\label{fig:halfkish}
\end{figure}

We  modify Kishino's knot and $ K_D $ to produce 
the diagram $ K_m $, shown in Figure \ref{fig:halfkish}, which shares two important characteristics with Kishino's knot. This knot is not detected by the bracket polynomial or fundamental group. We use this virtual knot diagram to construct an  infinite family (Family B) of non-trivial and non-classical virtual knot diagrams that are not detected by the fundamental group or bracket polynomial and are trivial as flats.

\begin{prop} \label{km0}
The non-trivial virtual knot diagram in Figure \ref{fig:halfkish} is not detected by the bracket polynomial or the 2-strand bracket polynomial, but is detected by the 3-strand bracket polynomial.
\end{prop}
\textbf{Proof:}
Note that:
\begin{alignat*}{2}
\langle K_m \rangle &= 1 & \qquad \langle K_m \rangle_2 &= d
\end{alignat*}
however
% Recompute and verify
\begin{gather*}
\langle K_m \rangle_3 =  -527 - A^{-60} + A^{-56} + 2A^{-54} +
5A^{-52} + 10A^{-50} \\
 +     21A^{-48} + 25A^{-46} + 26A^{-44} + 21A^{-42} + 8A^{-40} - 19A^{-38}\\
 -     69A^{-36} - 115A^{-34} - 155A^{-32} - 175A^{-30} - 172A^{-28} \\
 -     127A^{-26} - 23A^{-24} + 109A^{-22} + 244A^{-20} + 366A^{-18}\\
 + 
    440A^{-16} + 452A^{-14} + 372A^{-12} + 207A^{-10} + 23\/A^{-8} \\
- 199A^{-6} - 358A^{-4} - 486A^{-2} - 460A^2 - 364A^4 - 187A^6 \\
-  58A^8 + 96A^{10} + 208A^{12} + 250A^{14} + 269A^{16} \\
 +  222A^{18} + 200A^{20} + 131A^{22} + 71A^{24} + 21A^{26}  \\ - 
    15A^{28} - 29A^{30} - 54 A^{32} - 51A^{34} - 48A^{36} - 36A^{38} \\
 - 28A^{40} - 22 A^{42} - 10 A^{44} - 6A^{46} - 2A^{48} - A^{50} + A^{54}
\end{gather*}
Therefore, $ K_m $ is not equivalent to the unknot.\qed

\begin{figure}[htb] \epsfysize = 1 in
\centerline{\epsffile{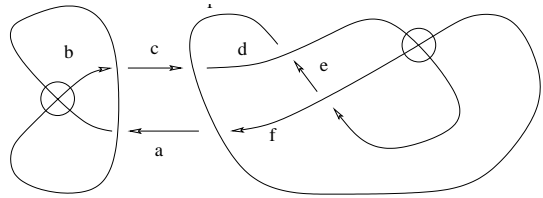}}
\caption{Fundamental Group of $ K_m $ }
\label{fig:modfund}
\end{figure}

\begin{prop}
The fundamental group of $ K_m $, shown in Figure \ref{fig:halfkish}, is $ \mathbb{Z}$.
\end{prop}
\textbf{Proof:}
To compute the fundamental group of  $ K_m $, we orient the knot and label each arc of the diagram. From each crossing, we obtain the one of the relations shown in Figure \ref{fig:modfund}.

The relations obtained from the left hand side of this diagram 
reduce to the equality: $ a= c $. This is the same relationship that would occur if the left hand side of this virtual diagram was replaced with an unknotted $ 1-1 $ tangle. The connected sum of an unknotted $ 1-1 $ tangle and the right hand side are equivalent to  unknot. 
 Hence, $  \pi_1 ( K_m ) $ reduces to  $ \mathbb{Z} $.  
\qed

\begin{figure}[htb] \epsfysize = 2.5 in
\centerline{\epsffile{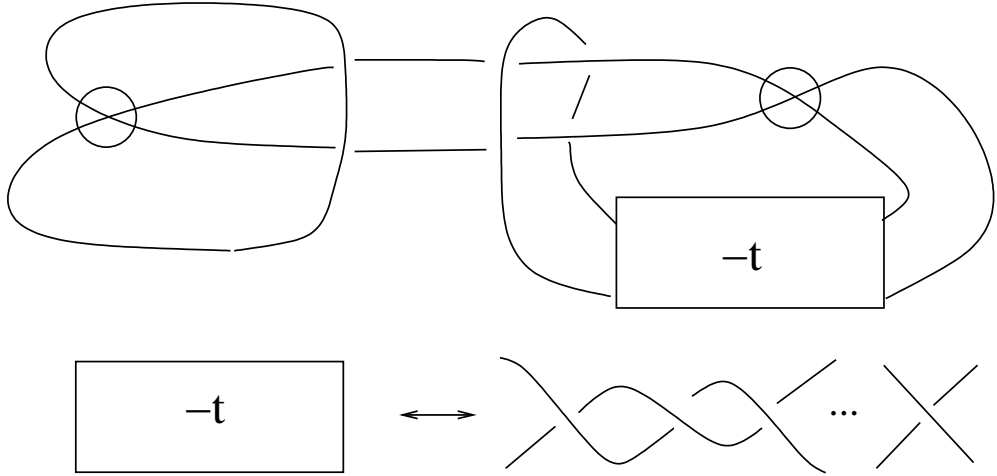}}
\caption{Schematic of Infinite Family B, $ K_m(t) $ }
\label{fig:mschem}
\end{figure}

From this diagram, we construct an infinite family of virtual knot diagrams. The members of this family are denoted $ K_m (t) $ 
where $ t $ represents the number of twists inserted into the virtual knot diagram $ K_m $. A schematic of this family is shown in Figure \ref{fig:mschem}. Note that $ K_m (0) $ denotes $ K_m $.

\begin{thm} For all $ t $, $ K_m(t) $, illustrated in the schematic diagram in Figure \ref{fig:mschem}, $ \langle K_m(t) \rangle = -A^{3t} $ and $ \pi_1 (K_m(t)) =  \mathbb{Z} $.
\end{thm}

\begin{figure}[htb] \epsfysize = 1.5 in
\centerline{\epsffile{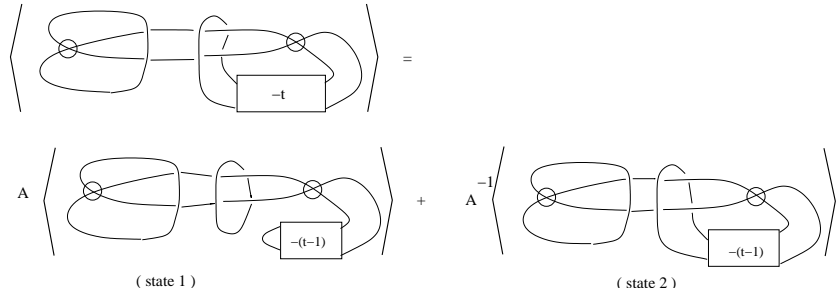}}
\caption{Expansion of $ K_m(t) $ }
\label{fig:mschemexpan}
\end{figure}

\textbf{Proof:} We compute $ \langle K_m (t) \rangle $. Recall that $ \langle K_m (0) \rangle = 1 $ as shown in Theorem \ref{km0}.
We expand the diagram $ K_m (t) $ using the skein relation as shown in Figure 
\ref{fig:mschemexpan}.

After a sequence of Reidemeister moves and virtual Reidemeister moves, state 1 is equivalent to two unknotted, unlinked components. One of these components has $ -(t-1) $ inserted twists. 
We note that $ \langle \text{state a} \rangle = -A^{3(t-1)} d $. State 2 is equivalent to $ K_m(t-1) $.
Hence
\begin{equation}
\langle K_m(t) \rangle = A \langle state 1 \rangle + A^{-1} \langle K_{m} (t-1) \rangle
\end{equation}
If we assume that $ \langle K_m (t-1) = -A^{3(t-1)} $ then 
\begin{equation*}
\langle K_m (t) \rangle = A ( -A)^{3(t-1)} d + A^{-1} (-A)^{3(t-1)}
\end{equation*}
This reduces to 
\begin{equation*}
\langle K_m(t) = (-A)^{3(t-1)}(A d + A^{-1}) = (-A)^{3t}
\end{equation*}
Since $ \langle K_m (0) \rangle =1 $, then 
$ \langle K_m (t) \rangle = (-A)^{3t} $ for all $ t $.

We show that $ \pi_1 (K_m (t) ) = \mathbb{Z} $ 
for all $ t $.
We consider the left hand side of the diagram $ K_m (t) $.
\begin{figure}[htb] \epsfysize = 1 in
\centerline{\epsffile{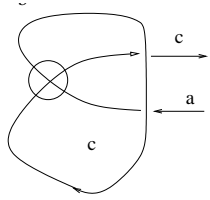}}
\caption{Left hand side of $ K_m(t) $ }
\label{fig:leftkm}
\end{figure}
Label the arcs as shown in Figure \ref{fig:leftkm}. The relations computed from this diagram reduce to:
\begin{gather*}
  a = b \text{ and }
   b  = c
\end{gather*}
Hence, we determine that  $ a= c $. 
The fundamental group of the left hand side is equivalent to the fundamental group of an unknotted $ 1-1 $ tangle. Taking the connected sum of the right hand side of the diagram $ K_m (t) $ and unknotted $ 1-1 $ tangle, we obtain a virtual knot diagram that is equivalent to the unknot. Thus, the fundamental group of $ K_m (t) $ is  $ \mathbb{Z} $.
\qed

\begin{thm}The virtual knot diagrams represented by the schematic diagram shown in Figure \ref{fig:mschem} are non-trivial and non-classical.
\end{thm}
\textbf{Proof:}

\begin{figure}[htb] \epsfysize = 1 in
\centerline{\epsffile{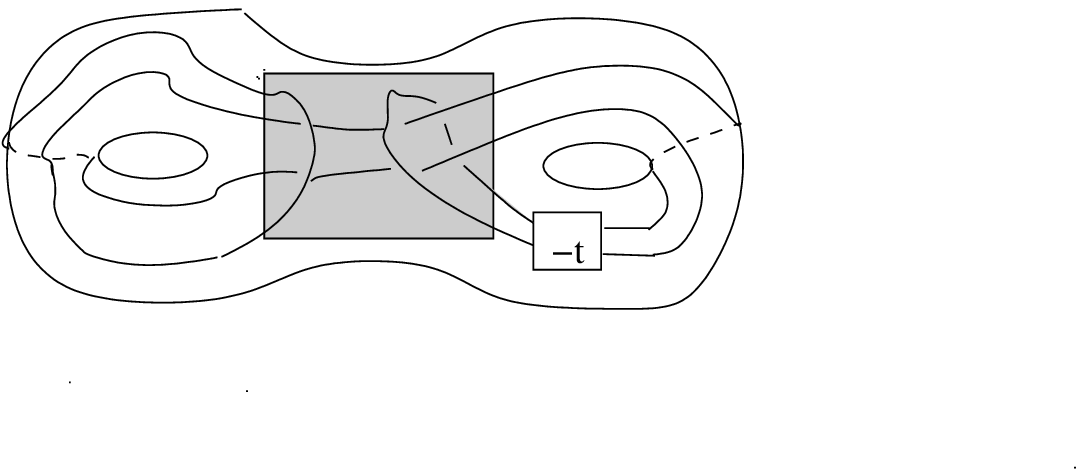}}
\caption{Schematic Representation of $ K_m (t)$}
\label{fig:kmsurf}
\end{figure}

In Figure \ref{fig:kmsurf}, we illustrate a schematic representation of $ K_m (t) $. We compute $\langle (F, K_m (0)) \rangle $.

\begin{figure}[htb] \epsfysize = 1 in
\centerline{\epsffile{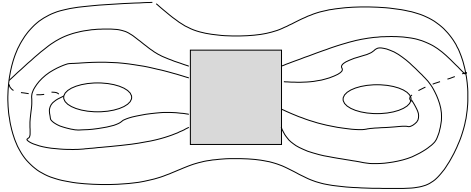}}
\caption{Expansion of $ K_m $}
\label{fig:kmexpand}
\end{figure}

Applying the skein relation to
expand the classical crossings in the shaded tangle box shown in Figure \ref{fig:kmsurf} results in 13 surface-state pairs with non-zero coefficients. 
These surface-state pairs are obtained by inserting
the possible expansions of $4-4 $ tangles from
Figure \ref{fig:expan4} into the shaded tangle box in Figure \ref{fig:kmexpand}.

Now:
\begin{gather*} 
\langle (F,K_m (0) ) \rangle =
A^2 \langle (F, s_2) \rangle + 
( - A^{-6} + A^2 ) \langle (F, s_3 ) \rangle 
+( - A^{-2}+A^2) \langle (F, s_4) \rangle \\
 +(- A^{-4}+ 1) \langle (F, s_5 ) \rangle
+( - A^{-4} + A^{4}) \langle (F, s_6) \rangle
+ ( -A^{-8} + 1) \langle (F, s_7) \rangle   \\    
+ A^4 \langle (F, s_8) \rangle  
+(- A^{-4} + 1 ) \langle (F, s_9) \rangle 
+(- A^{-2} + A^2 ) \langle (F, s_{10} ) \rangle  \\
+(- A^{-6} + A^2) \langle (F, s_{11} ) \rangle
+ (A^{-4} + 1 +A^4) \langle (F, s_{12} ) \rangle 
+ (-A^{-6}) +A^{-2}) \langle ( F, s_{13} ) \rangle  \\
+  A^6 \langle (F, s_{14} ) \rangle
\end{gather*}
The states $ (F, s_3) $, $(F, s_4)$, and $ (F, s_5) $ have non-zero coefficients,
indicating that the minimal genus of $ K_m (0) $ is two.
To compute the surface bracket polynomial of $ K_m (t) $ for $ t \geq 1 $, we note that we may apply the same expansion used in the proof of Theorem \ref{bkd}. 
We obtain two surfaces with embedded links.
The first surface-link pair admits a cancellation curve along the 
meridian of the right hand torus. The surface-state pairs obtained from applying the surface bracket polynomial to this surface-link pair do not include $ (F, s_3 ) $, $(F,s_4 ) $, or $ ( F, s_5) $.
The other surface-link pair is a representation of $ K_m (t-1) $.
Applying the skein relation to $ (F, K_m (t-1) ) $ results in the following equation:
\begin{gather*}
\langle (F, K_m (t) ) \rangle =
(-A)^{-4t} ( ( - A^{-6} + A^2 ) \langle (F, s_3 ) \rangle 
+( - A^{-2}+A^2) \langle (F, s_4) \rangle \\
 +(- A^{-4}+ 1) \langle (F, s_5 ) \rangle ) + X
\end{gather*}
where $ X $ represents all other possible surface state pairs
such as $ (F, s_7 ) $ and their coefficients.
where $ X $ denotes the other surface-state pairs and their coefficients. Hence, the virtual genus of $ K_m (t) $ is two. Therefore $ K_m (t) $ is non-classical and non-trivial for all $ t $.
\qed

\begin{rem}
With some modification, we can produce a related family that is not detected by the 2-strand bracket polynomial for small $ t $. We illustrate this family, denoted $ K_m'(t) $ in 
Figure \ref{fig:modfamB}.
\begin{figure}[htb] \epsfysize = 1 in
\centerline{\epsffile{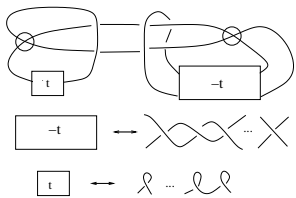}}
\caption{Family $ K_m'(t) $}
\label{fig:modfamB}
\end{figure}
\end{rem}

\begin{conj}We conjecture that $ K_m'(t) $ is a 
non-classical virtual knot diagram and is not detected by the fundamental group or the N-strand bracket 
polynomial unless $ N \geq 3 $.
\end{conj}
\section{Conclusion}

These new virtual diagrams provide a new benchmark in assessing the effectiveness and computability of virtual knot diagram invariants. 
We hope to consider the following questions:
\begin{itemize}
\item Do there exist virtual knot diagrams that are not detected by the 3-strand bracket polynomial?
\item Do there exist any tangles not detected by the N-strand bracket polynomial for 
a $ N \leq M $, M fixed ? 
\item Is there a geometric reason why the 3-strand bracket polynomial detects these new examples?
\end{itemize} 
We remark that the second question has been partially answered for the case when $ M = 2,3 $. 
In \cite{kvirt}, Kauffman determined that a single crossing, flanked by 
two virtual crossings was not detected by the bracket polynomial.

In response to the third question, it may be possible that the 3-strand 
bracket polynomial in some sense detects the minimum genus of these non-classical diagrams \cite{dk1}, \cite{kamada}.

We conclude this paper with the following conjecture:
\begin{conj}The 3-strand bracket polynomial detects all non-trivial virtual 
knot diagrams.
\end{conj}

\end{document}